\documentclass[12pt]{article}
\usepackage{amsmath,amssymb,latexsym}

\usepackage{color}

\definecolor{xxl}{rgb}{0.2,0.08,0.4}

\newtheorem{defn}{Definition}[section]
\newtheorem{lemma}[defn]{Lemma}
\newtheorem{ex}[defn]{Example}
\newtheorem{thm}[defn]{Theorem}
\newtheorem{prop}[defn]{Proposition}
\newtheorem{cor}[defn]{Corollary}
\newtheorem{rem}[defn]{Remark}

\newcommand{\h}{{\cal H}}
\newcommand{\mn}{\mathbb N}
\newcommand{\mr}{\mathbb R}

\def\range{{\cal R}}
\def\H{{\cal H}}
\def\h{{\cal H}}
\def\bp{\noindent{\bf Proof: \ }}
\def\ep{\noindent{$\Box$}}
\def\<{\langle}
\def\>{\rangle}

\newcommand{\seq[1]}{(#1_n)}

\def\newin {\,\kern-0.4em\in\kern-0.15em}
\def\newsubset {\kern-0.2em\subset\kern-0.2em}

\def\v{\vspace{.1in}}
\def\normsn{$\|\!\cdot\!\|$-$SN$ }

\title{Unconditional Convergence and \\
Invertibility of Multipliers}

{\small 
\author{
 D.\,T. Stoeva$^{a), b)}$\footnote{ std73std@yahoo.com} \ and P. Balazs$^{b)}$\footnote{ peter.balazs@oeaw.ac.at} \\ 
$^{a)}$ Department of Mathematics, \\ University of Architecture, Civil Engineering and Geodesy,\\
Blvd Hristo Smirnenski 1, 1046 Sofia, Bulgaria\\
$^{b)}$ Acoustics Research Institute, \\
Wohllebengasse 12-14, Vienna A-1040, Austria \\
}
}

\begin{document}
\maketitle \pagestyle{myheadings} 
{\footnotetext[1]{This work was supported by  the WWTF project MULAC ('Frame Multipliers: Theory and Application in Acoustics; MA07-025)} }

\parindent0pt
\parskip1ex

\begin{abstract}
In the present paper the unconditional convergence and the invertibility of multipliers is investigated. Multipliers are operators created by (frame-like) analysis, multiplication by a fixed symbol, and resynthesis. Sufficient and/or necessary conditions for unconditional convergence and invertibility are determined depending on the properties of the analysis and synthesis sequences, as well as the symbol. Examples which show that the given assertions cover different classes of multipliers are given.
If a multiplier is invertible, a formula for the inverse operator is determined.
The case when one of the sequences is a Riesz basis is completely characterized. 
\end{abstract}



\section{Introduction}

In modern life, applications of signal processing can be found in numerous technical items, for example in wireless communication or medical imaging. In these applications, \lq{\em time-invariant filters}\rq, i.e. convolution operators, are used very often. Such operators can be called {\em Fourier multipliers} \cite{fena06,bgrgrok05}. In the last decade {\em time-variant filters} have found more and more applications. A particular way to implement such filters are Gabor multipliers \cite{feinow1,benepfand07}, also known as Gabor filters \cite{hlawatgabfilt1}. Such operators find application in psychoacoustics \cite{xxllabmask1}, computational auditory scene analysis \cite{wanbro06}, virtual acoustics \cite{majxxl1}, and seismic data analysis \cite{gigrheillama05}. 
In \cite{xxlmult1} the concept of Bessel multipliers, i.e. operators of the form
 $$M f=\sum_k m_k \<f,\phi_k\>\psi_k, \ \forall f\in\h,$$ 
with $\seq[\phi]$ and $\seq[\psi]$ being Bessel sequences, were introduced and investigated. 
Further, the similar concept for $p$-Bessel sequences is considered in \cite{rahxxlXX}.
For many applications, for example in sound morphing \cite{DepKronTor07}, to find the inverse of such operators is of interest. In this paper we investigate the invertibility of multipliers.

From a theoretical point of view, it is very natural to investigate Bessel and frame multipliers. 
In \cite{schatt1}, R. Schatten investigated such operators for orthonormal sequences and provided a detailed
study of ideals of compact operators using their singular
decomposition.  By the spectral theorem, every self-adjoint compact operator on a Hilbert space can be represented as a  multiplier using an orthonormal system.
Moreover, multipliers generalize the frame operators, as every frame operator $S$ for a frame $\seq[\phi]$ is the multiplier $M_{(1),(\phi_n),(\phi_n)}$.

Multipliers have application as time-variant filters \cite{DepKronTor07,xxllabmask1,hlawatgabfilt1} in acoustical signal processing. Therefore, it is interesting to determine their inverses. 
For example, if the reverberating system of a church can be modeled by a multiplier, including time-variant aspects like e.g. the presence of an audience, having the inverse would give a possibility to get the original signal from a recorded signal. 
Furthermore, if some operator can be well approximated by a multiplier, we could solve an operator equation (like in e.g. \cite{kreizxxl1}) by inverting the multiplier.  
From a frame theory point of view, the investigation of multipliers in the special case of the identity implies properties of dual systems. 

Some properties of the invertibility of multipliers are known. For a frame $\seq[\phi]$ and a positive (resp. negative) semi-normalized sequence $\seq[m]$, the multiplier $M_{\seq[m],\seq[\phi],\seq[\phi]}$ is the frame operator $S$ (resp. $-S$) for the frame $(\sqrt{| m_n |} \,\phi_n)$ and thus, $M_{\seq[m],\seq[\phi],\seq[\phi]}$ is invertible \cite{xxljpa1}. When $\seq[\phi]$ and $\seq[\psi]$ are Riesz bases and $\seq[m]$ is semi-normalized, then $M_{\seq[m],\seq[\phi],\seq[\psi]}$ is invertible and $M_{\seq[m],\seq[\phi],\seq[\psi]}^{-1}=M_{(\frac{1}{m_n}),\seq[\widetilde{\psi}],\seq[\widetilde{\phi}]}$, where $\seq[\widetilde{\phi}]$ and $\seq[\widetilde{\psi}]$ denote the canonical duals of $\seq[\phi]$ and $\seq[\psi]$, respectively, see \cite{xxlmult1}.
 If $\seq[\phi^d]$ is a dual frame of the frame $\seq[\phi]$, then $M_{(1),\seq[\phi],\seq[\phi^d]}$ is the identity operator and therefore, invertible. 
If $m\in c_0$, and both $\seq[\phi]$ and $\seq[\psi]$ are Bessel sequences, then the multiplier $M_{\seq[m],\seq[\phi],\seq[\psi]}$ is never invertible on an infinite dimensional Hilbert space, because it is a compact operator \cite{xxlmult1}.

In the present paper, we investigate the invertibility of multipliers in more details. 
In Section \ref{sec:prel0}, we specify the notation and state the needed results for the main part of the paper. In Section \ref{uncconv}, the unconditional convergence of multipliers is considered; sufficient and/or necessary conditions 
are determined.  Section 4 concerns the question of the invertibility of multipliers $M_{(m_n),(\phi_n),(\psi_n)}$. Different cases for $(\phi_n)$ and $(\psi_n)$ are considered - non-Bessel, Bessel sequences, overcomplete frames, and Riesz bases. Sufficient and/or necessary conditions for the invertibility of $M_{(m_n),(\phi_n),(\psi_n)}$ are given. 
If the multipliers are invertible, formulas for $M_{(m_n),(\phi_n),(\psi_n)}^{-1}$ are determined. 
It is planned to use those in numerical algorithms in the future. 
The last section of the paper contains examples showing the sharpness of the bounds in Propositions \ref{p1} - 
\ref{p3}, as well as showing that Propositions \ref{mpos}, \ref{p4} and \ref{p3} are independent of each other.

For certain cases of multipliers we provide examples and counter-examples. For some of them we refer to the paper \cite{BStable09}.
Our goal there was to characterize a complete set of conditions for the invertibility and unconditional convergence of multipliers, beyond the set of examples needed for the current paper. 
The complete collection of these examples can be found in  \cite{BStable09}.

\section{Notation and preliminary results } \label{sec:prel0}

Throughout the paper 
$\h$ denotes a Hilbert space and $\seq[e]$ denotes an orthonormal basis of $\h$. 
The notion {\it operator} is used for linear mappings. The range of an operator $G$ is denoted by $\range(G)$. 
The identity operator on $\h$ is denoted by $I_\h$. The operator $G:\h\to \h$ is called {\it invertible on $\h$} if 
 there exists a bounded operator $G^{-1}:\h\to \h$  such that $GG^{-1}=G^{-1}G=I_\h$ (and therefore $G$ is bounded). 
Throughout the paper, we work with a fixed infinite, but countable index set $J$, and, without loss of generality, $\mn$ is used as an index set, also implicitly.

The notation $\Phi$ (resp. $\Psi$) is used to denote the sequence $\seq[\phi]$ (resp. $\seq[\psi]$) with elements from $\h$; $\Phi-\Psi$ denotes the sequence $(\phi_n-\psi_n)$; $m$ denotes a 
 complex scalar sequence $\seq[m]$,
 and $\overline{m}$ denotes the sequence of the complex conjugates of $m_n$; $m\Phi$ denotes the sequence $(m_n\phi_n)$.  Recall that $m$ is called {\it semi-normalized} (in short, $SN$) if there exist constants $a,b$ such that $0<a\leq |m_n|\leq b<\infty$, $\forall n$. A series $\sum \phi_n$ is called {\it unconditionally convergent} if 
$\sum \phi_{\sigma(n)}$   converges for every permutation $\sigma(n)$ of $\mn$.
   
Recall that $\Phi$ is called a {\it Bessel sequence} (in short, {\it Bessel}) {\it for $\h$ with bound $B_\Phi$} if $B_\Phi>0$ 
and $\sum |\<h,\phi_n\>|^2 \leq B_\Phi\|h\|^2$ for every $h\in\h$. A Bessel sequence $\Phi$ with bound $B_\Phi$  is called a {\it frame for $\h$ with bounds $A_\Phi, B_\Phi$}, if $A_\Phi>0$ and $A_\Phi\|h\|^2\leq  \sum |\<h,\phi_n\>|^2 $ for every $h\in\h$; $A_\Phi^{opt}$ and $B_\Phi^{opt}$ denote the optimal frame bounds for $\Phi$. 
The sequence $\Phi$ is called a Riesz basis for $\h$ with bounds $A_\Phi,B_\Phi$,  
if $A_\Phi>0$ and 
$A_\Phi \sum |c_n|^2\leq \|\sum c_n \phi_n\|^2\leq B_\Phi \sum |c_n|^2$,
$\forall \seq[c]\in\ell^2$.
Every Riesz basis for $\h$ with bounds $A,B$ is a frame for $\h$ with bounds $A,B$. 
For standard references for frame theory and related topics see \cite{Casaz1,ole1,he98-1}.

For a given Bessel sequence $\Phi$, the mapping $U_\Phi:\h\to \ell^2$ given by $U_\Phi f =(\<f,\phi_n\>)$ is called {\it the analysis operator for $\Phi$} and the mapping $T_\Phi$ given by $T_\Phi \seq[c]= \sum c_n \phi_n$   is called {\it the synthesis operator for $\Phi$}. Let $\Phi$ be a frame for $\h$. The operator $S_\Phi:\h\to\h$ given by $S_\Phi h=\sum\<h,\phi_n\>\phi_n$ is called the {\it frame operator for $\Phi$} and fulfills  
$A_\Phi\|f\|\leq \|Sf\|\leq B_\Phi\|f\|$, $\forall f\in\h$.
The sequence $(\phi_n^d)$ is called a {\it dual frame of $\Phi$} if it is a frame for $\h$ and $h=\sum \<h,\phi_n^d\>\phi_n=\sum \<h,\phi_n\>\phi_n^d$, $\forall h\in\h$. The sequence $\widetilde{\Phi}=(S^{-1}_\Phi \phi_n)$ is a dual frame of $\Phi$, called {\it the canonical dual of $\Phi$}. 
Recall the following results:

\begin{prop} \label{Cprop} For a sequence $\Phi$, the following statements hold.
\
\begin{itemize}
\item[{\rm (i)}] {\rm \cite[p.75]{he98-1}}  $\Phi$ is Bessel for $\h$ if and only if $(\<f,\phi_n\>)\in\ell^2$ for every $f\in\h$.
\item[{\rm (ii)}] {\rm \cite[p.52]{ole1}} 
$\Phi$ is Bessel for $\h$ if and only if the operator $T_\Phi : \seq[c]\to \sum c_n \phi_n$ is well defined from $\ell^2$ into $\h$; in this case $\|T_\Phi\|= \sqrt{B_\Phi^{opt}}$.
\item[{\rm (iii)}] {\rm \cite[Prop. 12.10]{he98-1}}
If $\Phi$ is a frame for $\h$ with bounds $A_\Phi,B_\Phi$, then 
$S_\Phi$ is invertible on $\H$
 and $\frac{1}{B_\Phi} \|f\| \leq \|S_\Phi^{-1} f\|\leq\frac{1}{A_\Phi} \|f\|$, $\forall f\newin\h$. 
\end{itemize}
\end{prop}

If $\sup_n \|\phi_n\| <\infty$ (resp. $\inf_n \|\phi_n\| >0)$, the sequence $\Phi$ will be called {\it norm-bounded above}, in short $NBA$  (resp. {\it norm-bounded below}, in short $NBB$). 
If $(\|\phi_n\|)$ is semi-normalized, then $\Phi$ is called {\it $\|\!\cdot\!\|$-semi-normalized} (in short, $\|\!\cdot\!\|$-SN). 
Recall that if $\Phi$ is a Bessel sequence for $\h$ with bound $B_\Phi$, then $\|\phi_n\|\leq \sqrt{B_\Phi}$, $\forall n$, and clearly, $\Phi$ is not needed to be $NBB$. 
Note that even a frame $\Phi$ is not needed to be $NBB$. Take, for example, the frame $(\frac{1}{2} e_1, e_2,
\frac{1}{2^2} e_1, e_3, \frac{1}{2^3} e_1, e_4, \ldots)$.
 Typical examples for \normsn frames are Gabor and wavelet frames, \cite{gr01}.  Any Riesz basis $\Phi$ is  \normsn, because $\sqrt{A_\Phi}\leq \|\phi_n\|\leq \sqrt{B_\Phi}, \forall n$.
 A frame which is \normsn does not need to be a Riesz basis; consider, for example, the sequence $(e_1,e_1,e_2,e_3,\ldots)$. 
We need the following result concerning Riesz bases:

\begin{prop} {\rm \cite{stoevxxl09}} \label{ovfr} \
The sequence $m\Phi$ can be a Riesz basis for $\h$ only in the following cases:

\hspace{.15in} $\bullet$ $\Phi$ is a Riesz basis for $\h$ and $m$ is $SN$;

\hspace{.15in} $\bullet$ $\Phi$ is non-$NBB$ Bessel for $\h$ which is not a frame for $\h$ and $m$ is $NBB$ but not in $\ell^\infty$;

\hspace{.15in} $\bullet$ $\Phi$ is non-$NBA$ non-Bessel for $\h$ and $m$ is non-$NBB$ with $m_n\neq 0$, $\forall n$.

\end{prop}

  Note that the paper \cite{stoevxxl09} concerns real sequences $m$, but the above statement holds also for complex sequences $m$.
  Further, it is easy to observe the following relationship between sequences $m\Phi$ and $\overline{m}\,\Phi$:

\begin{prop} \label{mmbar}
The sequence $m\Phi$ is Bessel, frame, Riesz basis or satisfies the lower frame condition for $\h$, if and only if $\overline{m}\,\Phi$ is Bessel, frame, Riesz basis or satisfies the lower frame condition for $\h$, respectively, with the same bounds. 
\end{prop}

\subsection*{Multipliers}

For any $\Phi$, $\Psi$ and any $m$ (called {\it weight} or {\it symbol}), the operator $M_{m,\Phi,\Psi}$, given by $$M_{m,\Phi, \Psi} f=\sum m_n \<f,\psi_n\> \phi_n, \ f\newin\h,$$ is called a {\it multiplier} \cite{xxlmult1}. When $\Phi$ and $\Psi$ are Bessel sequences (resp. frames), $M_{m,\Phi, \Psi}$ is called a {\it Bessel} (resp. {\it frame}) {\it multiplier}. Depending on $m,\Phi$, and $\Psi$, the corresponding multiplier might not be well defined, i.e. might not converge for some $f\in\h$. 
The following assertion gives a sufficient condition for the well-definedness of multipliers. 

\begin{prop} {\rm \cite{xxlmult1}} \label{bmh}
Let $m\newin\ell^\infty$ and let $\Phi$, $\Psi$ be Bessel sequences for $\h$ with bounds $B_\Phi,B_\Psi$, respectively. Then the multiplier $M_{m,\Phi,\Psi}$
is well defined from $\h$ into $\h$ and it is bounded with $\|M_{m,\Phi,\Psi} \|\leq \sqrt{B_\Phi\,B_\Psi}\,\|m\|_\infty$. Furthermore, the series $\sum m_n \<f,\psi_n\> \phi_n$ converges unconditionally for every $f\newin\h $.
\end{prop}

Note that the above assertion gives only a sufficient condition. Multipliers can be unconditionally convergent and bounded even in cases when $m\notin\ell^\infty$ or at least one of the sequences is not Bessel. For example, consider $M_{(n^2),(\frac{1}{n}e_n),(\frac{1}{n}e_n)}=I_\h$ and $M_{(1),(\frac{1}{n}e_n),(ne_n)}=I_\h$.
The question for unconditional convergence of multipliers is investigated deeper in Section \ref{uncconv}. 
We will use the following known results:

\begin{prop}  \label{hprop} For a sequence $\Phi$, the following statements hold.
\begin{itemize} 
\item[{\rm (i)}] {\rm \cite[p.32]{he98-1}} If $\sum \phi_n$ converges unconditionally, then
$\sum \|\phi_n\|^2<\infty$. 

\item[{\rm (ii)}] {\rm \cite{Orl, Pettis, he98-1}} The following conditions are equivalent.

$\bullet$ $\sum_n \phi_n$ converges unconditionally.

$\bullet$ Every subseries $\sum_k \phi_{n_k}$ converges.

$\bullet$ Every subseries $\sum_k \phi_{n_k}$ converges weakly.

$\bullet$ $\sum_n \lambda_n \phi_n$ converges for every bounded sequence of scalars $(\lambda_n)$.

\item[{\rm (iii)}] {\rm \cite[p.92]{he98-1}} If $\Phi$ is a Riesz basis for $\h$, then 
$\sum c_n \phi_n$ converges unconditionally  if and only if $\sum c_n \phi_n$ converges.

\item[{\rm (iv)}] If $\Phi$ is a $NBB$ Bessel sequence for $\h$, then $\sum c_n \phi_n$ converges unconditionally if and only if $\seq[c]\in\ell^2$.  

\end{itemize}
\end{prop}
If $\Phi$ is a $NBB$ frame for $\h$, the conclusion of Proposition \ref{hprop}(iv) is proved in \cite[Prop. 12.17]{he98-1}. The proof in \cite{he98-1} uses only validity of the upper frame condition, so the property is shown for Bessel sequences.

Concerning Proposition \ref{hprop}(iv), note that if the condition \lq\lq norm-bounded below\rq\rq \ is omitted, then the conclusion does not hold in general, because $\sum c_n \phi_n$ might converge unconditionally for some $\seq[c]\notin\ell^\infty$, see \cite[Ex. 12.16]{he98-1}.

\v The main aim of our paper is to investigate the invertibility of multipliers. We will use the following criterion 
for the invertibility of operators:

\begin{prop}  \label{ginvnew}
Let $F:\h\to\h$ be invertible on $\h$. Suppose that $G:\h\to\h$ is a bounded operator and $\|Gh - Fh\|\leq \nu \|h\|$, $\forall h\in \h$, where $\nu\in [0, \frac{1}{\|F^{-1}\|})$. 
Then 
\begin{itemize}
\item[{\rm (i)}] $G$ is invertible on $\h$ and 
$G^{-1}=\sum_{k=0}^\infty [F^{-1}(F-G)]^kF^{-1};$
\item[{\rm (ii)}] 
$\frac{1}{1+\nu \|F^{-1}\|} \,\frac{1}{\|F\|}\,\|h\| \leq \|G^{-1}h\|\leq \frac{1}{\left(\frac{1}{\|F^{-1}\|} \, - \, \nu\right)} \, \|h\|, \ \forall h\in \h.$

\end{itemize}
\end{prop}
\bp
(i) is proved in {\rm \cite[p.70]{gohberg1}}. 
For (ii), observe that 
$\|Gh-Fh\|\leq \nu \,\|F^{-1}\|\, \|Fh\|$  and apply  \cite[Theorem 1]{cach97} with $\lambda_1= \nu \|F^{-1}\|<1$ and $\lambda_2=0$.
\ep

\section{Unconditional convergence of multipliers}\label{uncconv}

First, observe the following easy consequence of Banach-Steinhaus Theorem:
\begin{lemma}
If $M_{m,\Phi,\Psi}$ is well-defined on $\h$, then it is bounded.
\end{lemma}

Further, we continue with a stronger notion of convergence.
The multiplier $M_{m,\Phi,\Psi}$ is called {\it unconditionally convergent on $\h$} 
if $\sum m_n\<f,\psi_n\>\phi_n$ converges unconditionally for every $f\in\h$.
First, we prove that the unconditional convergence of $M_{m,\Phi,\Psi}$ on $\h$ is equivalent to the unconditional convergence of $M_{m,\Psi,\Phi}$ on $\h$.

\begin{lemma} For a multiplier $M_{m,\Phi,\Psi}$ the following statements hold. 
\begin{itemize}
\item[{\rm (i)}] $M_{m,\Phi,\Psi}$ is unconditionally convergent on $\h$ if and only if $M_{\overline{m},\Psi,\Phi}$ is unconditionally convergent on $\h$.
\item[{\rm (ii)}] For any $f\in\h$, $M_{\overline{m},\Psi,\Phi}f$ converges unconditionally if and only if $M_{m,\Psi,\Phi}f$ converges unconditionally.
\end{itemize}
\end{lemma}
\bp (i) Let $M_{m,\Phi,\Psi}$ be unconditionally convergent on $\h$.  
By Proposition \ref{hprop}(ii), every subseries $\sum_k m_{n_k} \<f,\psi_{n_k}\> \phi_{n_k}$ converges for every $f\in \h$, which implies that every subseries 
$\sum_k \overline{m}_{n_k} \<g,\phi_{n_k}\> \psi_{n_k}$ converges weakly for every $g\in \h$. Now Proposition \ref{hprop}(ii) implies that $\sum_n \overline{m}_n \<g,\phi_n\> \psi_n$ converges unconditionally for every $g\in\h$.

(ii) Let $f\in\h$ and let $M_{\overline{m},\Psi,\Phi}f$ be unconditionally convergent. 
Then every subseries $\sum_k \overline{m}_{n_k} \<f,\phi_{n_k}\> \psi_{n_k}$ converges unconditionally. 
Consider the sequence $(\lambda_n)$ given by $\lambda_n=\frac{m_n}{\overline{m}_n}$ if $m_n\neq 0$ and $\lambda_n=0$ if $m_n=0$.
Applying Proposition \ref{hprop}(ii) with the bounded sequences $(\lambda_{n_k})_k$, it follows that every subseries $\sum_k m_{n_k} \<f,\phi_{n_k}\> \psi_{n_k}$ 
 converges. Now apply again Proposition \ref{hprop}(ii). 
 \ep

As a consequence, the following statement holds.
\begin{prop} \label{uncnew} For any sequences $m$, $\Phi$, and $\Psi$, the multiplier
 $M_{m,\Phi,\Psi}$ is unconditionally convergent on $\h$ if and only if $M_{m,\Psi,\Phi}$ is unconditionally convergent on $\h$.
\end{prop}

\begin{rem} 
Note that conditional convergence of $M_{m,\Phi,\Psi}$ is not equivalent to conditional convergence of $M_{m,\Psi,\Phi}$.
Consider for example the sequences 

\ \ \ $\Phi=(e_1, \ \, e_1, \ \ e_1, \ \  e_2,  \ \, e_2, \ \ e_2, \ \ e_3,  \  \, e_3, \ \  e_3,   \ldots)$,

\ \ \ $\Psi=(e_1, \ \, e_1,   - e_1, \ \,  e_2,  \ \, e_1, - e_1,  \ \ e_3,  \ \, e_1,     - e_1,   \ldots)$. 

\noindent Then $M_{m,\Psi,\Phi}=I_\h$ and $M_{m,\Phi,\Psi}$ is not well-defined. 
\end{rem}

As one can see in Proposition \ref{bmh}, Bessel multipliers are unconditionally convergent in case $m\in\ell^\infty$. Now we are interested in converse assertions.

\begin{lemma} \label{lem31} Let $M_{m,\Phi,\Psi}$ be unconditionally convergent on $\h$. 
\begin{itemize}
\item[{\rm (i)}] Then $(|m_n|\cdot \|\phi_n\| \cdot \psi_n)$ and $(|m_n|\cdot \|\psi_n\| \cdot \phi_n)$ are Bessel for $\h$. 
As a consequence, the sequence $(|m_n| \cdot \|\phi_n\| \cdot \|\psi_n\|)$ is bounded.
\item[{\rm (ii)}] If $\Phi$ (resp. $m\Phi$) is $NBB$, then $m\Psi$ (resp. $\Psi$) is a Bessel sequence for $\h$.
\item[{\rm (iii)}] If $\Psi$ (resp. $m\Psi$) is $NBB$, then $m\Phi$ (resp. $\Phi$) is a Bessel sequence for $\h$.
\item[{\rm (iv)}] 
If both $\Phi$ and $\Psi$ are $NBB$, then $m\in\ell^\infty$.
\item[{\rm (v)}] 
If $\Phi$, $\Psi$ and $m$ are $NBB$, then $m$ is $SN$ and both $\Phi$ and $\Psi$ are Bessel sequences for $\h$. 
\end{itemize}
\end{lemma}
\bp (i)
It follows from Proposition \ref{hprop}(i) that $( \<f,|m_n| \cdot \|\phi_n\|\cdot \psi_n\> )\in\ell^2$ for every $f\newin\h$. Proposition \ref{Cprop}(i) implies that $(|m_n| \cdot \|\phi_n\|\cdot\psi_n)$ is Bessel for $\h$.
Now use Proposition \ref{uncnew} and apply what is already proved to $M_{m,\Psi,\Phi}$.

(ii)-(iv) follow easily from (i); (v) follows from (ii)-(iv).
\ep

\begin{cor} \label{nbunc}
If $\Psi$ is non-Bessel for $\h$ and $m\Phi$ is NBB, then $M_{m,\Phi,\Psi}$ can not be unconditionally convergent on $\h$.
As a consequence (using Proposition \ref{hprop}(iii)), if $\Psi$ is non-Bessel for $\h$ and $m\Phi$ is a Riesz basis for $\h$, then $M_{m,\Phi,\Psi}$ is not well defined.
\end{cor}

\begin{rem} Note that if $\Psi$ is non-Bessel for $\h$ and  $m\Phi$ is $NBB$ non-Bessel for $\h$, then $M_{m,\Phi,\Psi}$ can be conditionally convergent and invertible on $\h$, although $M_{m,\Phi,\Psi}$ can not be unconditionally convergent on $\h$ by Corollary \ref{nbunc}.
 Consider, for example, the non-Bessel sequences

\ \ \ $\Phi=(e_1, e_2, e_2,  -e_2,   e_3, e_3, -e_3,   e_3, -e_3,     e_4, e_4, -e_4,    e_4, -e_4,  e_4, -e_4, \ldots)$,

\ \ \ $\Psi=(e_1, e_2, e_2, \ \ e_2, e_3, e_3, \ \ e_3, e_3, \ \ e_3, e_4, e_4, \ \ \, e_4, e_4, \ \ e_4, e_4, \ \ e_4, \ldots)$, 

\noindent for which $M_{(1),\Phi,\Psi}=M_{(1),\Psi,\Phi}=I_\h$. 
\end{rem}

\begin{rem} \label{nbnb} Note that if $\Psi$ is non-Bessel for $\h$ and $m\Phi$ is non-$NBB$, then $M_{m,\Phi,\Psi}$ can be unconditionally convergent and invertible on $\h$. Consider, for example, the non-Bessel sequences

\ \ \ $\Phi=(\frac{1}{2}\,e_1,   2\,e_2,    \frac{1}{2^2}\,e_1,  3\,e_3, \frac{1}{2^3}\,e_1, 4\,e_4, \ldots)$, 
$\Psi=(e_1, \, \frac{1}{2}\, e_2,  \, e_1, \, \frac{1}{3}\, e_3,  \,  e_1,  \, \frac{1}{4}\, e_4, \ldots)$. 

Then $M_{(1),\Phi,\Psi}=M_{(1),\Psi,\Phi}=I_\h$ with unconditional convergence on $\h$, because 
 $M_{(1),\Phi,\Psi}=M_{(1),(c_n\phi_n),(\frac{1}{c_n}\psi_n)}$ and $M_{(1),\Psi,\Phi}=M_{(1),(\frac{1}{c_n}\psi_n),(c_n\phi_n)}$, where
 $(c_n)=  (\sqrt{2},   \, \frac{1}{2},   \sqrt{2^2},  \, \frac{1}{3}, \sqrt{2^3}, \, \frac{1}{4}, \ldots)$,  and the sequences $(c_n\phi_n)$, $(\frac{1}{c_n}\psi_n)$ are Bessel for $\h$ (apply Prop. \ref{bmh}).
\end{rem}

In certain cases, one can determine conditions which are necessary and sufficient for the unconditional convergence of  multipliers:

\begin{prop} \label{lem32} For a multiplier $M_{m,\Phi,\Psi}$, the following statements hold.
\begin{itemize}
\item[{\rm (i)}] Let $\Phi$ be a $NBB$ Bessel sequence for $\h$. Then 

$M_{m,\Phi,\Psi}$ is unconditionally convergent on $\h$ $\Leftrightarrow$ 
$m\Psi$ is Bessel for $\h$.

\item[{\rm (ii)}] Let $\Phi$ be a Riesz basis for $\h$. Then 

$M_{m,\Phi,\Psi}$ is well defined on $\h$ $\Leftrightarrow$   $M_{m,\Phi,\Psi}$ is unconditionally convergent on $\h$  
$\Leftrightarrow$  $m\Psi$ is Bessel for $\h$  $\Leftrightarrow$
 $M_{m,\Psi,\Phi}$ is well defined on $\h$
 $\Leftrightarrow$ $M_{m,\Psi,\Phi}$ is unconditionally convergent on $\h$. 
\item[{\rm (iii)}] Let $\Phi$ be a Riesz basis for $\h$ and $\Psi$ be $NBB$. Then 

$M_{m,\Phi,\Psi}$ (or $M_{m,\Psi,\Phi}$) is well defined on $\h$ 
$\Rightarrow$ $m\in\ell^\infty$. The converse does not hold in general.

\item[{\rm (iv)}] 
If $\Phi$ and $\Psi$ are Riesz bases for $\h$, then 
$M_{m,\Phi,\Psi}$ is well defined on $\h$ if and only if $m\in\ell^\infty$. 
\end{itemize} 
\end{prop}

\bp (i) 
By Proposition \ref{hprop}(iv), $M_{m,\Phi,\Psi}$ is unconditionally convergent on $\h$ if and only if $(\<f,\overline{m}_n\psi_n\>)\in\ell^2$, $\forall f\in\h$, which by Propositions \ref{Cprop}(i) and \ref{mmbar} is equivalent to
$m\Psi$ being Bessel for $\h$. 

(ii) 
The first equivalence follows from Proposition \ref{hprop}(iii). The second equivalence follows from (i), because Riesz bases are $NBB$ Bessel sequences.

For the third equivalence, consider $M_{m,\Psi,\Phi} f=\sum \<f, \phi_n\> m_n \psi_n,$ $f\newin\h$. 
By Proposition \ref{Cprop}(ii), the sequence $m\Psi$ is Bessel for $\h$ if and only if $\sum c_n m_n\psi_n$ converges for every $\seq[c]\in\ell^2$.
This holds if and only if $\sum \<f, \phi_n\> m_n \psi_n$ converges for every $f\newin\h$, because $\ell^2=\{(\<f,\phi_n\>) : f\in \h\}$ as $\Phi$ is a Riesz basis for $\h$ \cite[Prop. 5.1.5]{gr01}. 

To complete the last equivalence, use Proposition \ref{uncnew}.

(iii)  Assume that $M_{m,\Phi,\Psi}$  is well defined, or equivalently, by (ii), that  $M_{m,\Psi,\Phi}$ is well defined.
Let $a_\Psi>0$ denote a lower bound for $(\|\psi_n\|)$.  By (ii),  $m\Psi$ is Bessel for $\h$. Then $a_\Psi |m_n|\leq \|m_n\psi_n\|\leq \sqrt{B_{m\Psi}}$, which implies that $m$ belongs to $\ell^\infty$. 
For the converse, consider the multiplier $M_{(\frac{1}{n}), (e_n),(n^2 e_n)}$, which is not well defined.

(iv) follows form (iii) and  Proposition \ref{bmh}.
\ep
\begin{cor} \label{corlem32}
If it is moreover assumed that $m$ is $NBB$ (resp. $SN$), then each of the equivalent assertions in Proposition \ref{lem32}(i),(ii) implies (resp. is equivalent to) $\Psi$ being Bessel for $\h$. 
\end{cor}

\begin{rem} \label{cex1} If $\Phi$ is Bessel for $\h$, which is non-$NBB$, then the conclusion of Proposition \ref{lem32}(i) might fail. Consider $\Phi=(\frac{1}{2}e_1,e_2,\frac{1}{2^2}e_1,e_3,\frac{1}{2^3}e_1,e_4,\ldots)$, which is Bessel for $\h$, and  $\Psi=(e_1,e_2,e_1,e_3,e_1,e_4,\ldots)$, which is non-Bessel for $\h$. Then $M_{(1),\Phi,\Psi}=M_{(1),\Psi,\Phi}=I_\h$   
with unconditional convergence on $\h$, because $M_{(1),\Phi,\Psi}=M_{(1),(c_n\phi_n),(\frac{1}{c_n}\psi_n)}$ and $M_{(1),\Psi,\Phi}=M_{(1),(d_n\psi_n),(c_n\phi_n)}$, where
 $(c_n)=(\sqrt{2},   \, 1,   \sqrt{2^2},  \, 1, \sqrt{2^3}, \, 1, \ldots)$, and the sequences $(c_n\phi_n)$, $(\frac{1}{c_n}\psi_n)$ are Bessel for $\h$ (apply Prop. \ref{bmh}). But $m\Psi=\Psi$ is not Bessel for $\h$.

\end{rem}

 \section{Invertibility of multipliers} \label{s4}

First note that having zero elements at \lq\lq appropriate places\rq\rq\ of $\Phi$, $\Psi$ and $m$, one can get any desired multiplier, for example, the invertible identity operator and the zero 
 operator:

\begin{ex} 
{\bf (ZERO-example)} \label{zeroex}
\

Let $\Phi=(*,0,*,0,*,0,\ldots)$ and $\Psi=(0,*,0,*,0,*,\ldots)$, where the stars denote arbitrary elements of the space (not necessarily the same).
For any weight $m$, we have $M_{m,\Phi,\Psi}=M_{m,\Psi,\Phi}=0$ and thus, both multipliers are non-invertible on $\h$.

\end{ex}
 Note that in this example the sequences $\Phi$ and $\Psi$ can be any kind (e.g.: non-Bessel, Bessel non-frame, overcomplete frames, etc.) except Riesz bases.

\begin{ex} {\bf (IDENTITY-example)} \label{idex}
Let $\Phi=(*,e_1,0,*,e_2,0,*,e_3,0,\ldots)$,
 $\Psi=(0,e_1,*,0,e_2,*,0,e_3,*,\ldots)$ and\, $m=(*,1,*,*,1,*,*,1,*,\ldots)$. Then $M_{m,\Phi,\Psi} = M_{m,\Psi,\Phi}=I_\h$ and thus, both multipliers are invertible on $\h$. 
\end{ex}

Since the zero elements in the sequences $m,\Phi,$ and $\Psi$, do not have an influence on the values of the corresponding multiplier $M_{m,\Phi,\Psi}$, 
instead of the initial index set $J$ we can consider a new index set $J_0=J\setminus \{i : m_i=0 \mbox{ or } \phi_n=0 \mbox{ or } \psi_n=0\}$%
\footnote{As an outlook it can be interesting to link such a construction with the excess of frames \cite{bacahela03,bacahela03-1}.}%
.
Note that if $J_0$ is empty or finite and $\h$ is infinite dimensional, then $M_{m,\Phi,\Psi}$ can not be surjective and thus, it can not be invertible on $\h$. That is why only infinite $J_0$ is of interest for the present paper.
Without loss of generality, from now on we consider only sequences $m,\Phi,$ and $\Psi$, which do not contain zero elements, and $\mn$ is the index set.

Observe that if $M_{m,\Phi,\Psi}$ is invertible on $\h$, then $\Phi$ must be complete in $\h$.

\subsection{Multipliers for non-Bessel sequences} 

First note that it is possible to have an invertible unconditionally convergent multiplier even in cases when both sequences $\Phi$ and $\Psi$ are non-Bessel. Consider the trivial example  $M_{(\frac{1}{n^2}),(ne_n),(ne_n)}=I_\h$. 
Moreover, this is possible even in cases with $m=(1)$; see the sequences in Remark \ref{nbnb}.

Having in mind Corollary \ref{nbunc}, in the cases when $\Psi$ is non-Bessel for $\h$ considering unconditionally convergent multipliers $M_{m,\Phi,\Psi}$ is only possible if $m\Phi$ is non-$NBB$ non-Bessel or non-$NBB$ Bessel (in particular, could be a frame, but not a Riesz basis) - for any of these cases invertible and non-invertible multipliers exist, see \cite{BStable09}.

\subsection{Sufficient and/or necessary conditions for invertibility of Bessel multipliers}

\v If the multiplier $M_{(1),\Phi,\Psi}$ is invertible and one of the sequences $\Psi$ and $\Phi$ is Bessel, then the other one does not need to be Bessel. For example, consider the sequences in Remark \ref{cex1}.
Below we observe that if one of the sequences is Bessel, invertibility of $M_{(1),\Phi,\Psi}$ implies that the other one must satisfy the lower frame condition.

\begin{thm} \label{c1}
Let $M_{m,\Phi,\Psi}$ be invertible on $\h$. 

\begin{itemize}
\item[{\rm (i)}] If $\Psi$ (resp. $\Phi$) is a Bessel sequence for $\h$ with bound $B$, then $m\Phi$ (resp. $m\Psi$)  satisfies the lower frame condition for $\h$ with bound $\frac{1}{B\, \|M_{m,\Phi,\Psi}^{-1}\|^2}$.

\item[{\rm (ii)}] If $\Psi$ (resp. $\Phi$) and $m\Phi$ (resp. $m\Psi$) are Bessel sequences for $\h$, then they are frames for $\h$.
\item[{\rm (iii)}] If $\Psi$ (resp. $\Phi$) is a Bessel sequence for $\h$ and $m\in\ell^\infty$, then $\Phi$ (resp. $\Psi$) satisfies the lower frame condition for $\h$.
\item[{\rm (iv)}]
If $\Psi$ and $\Phi$ are Bessel sequences for $\h$ and $m\in\ell^\infty$, then $\Psi$ and $\Phi$ are frames for $\h$; $m\Phi$ and $m\Psi$ are also frames for $\h$. 
\end{itemize}
\end{thm}
\bp (i)  For brevity, the multiplier $M_{m,\Phi,\Psi}$ will be denoted by $M$. The proof will be done in two steps.

First step: $m=(1)$.

 Assume that $\Psi$ is a Bessel sequence for $\h$ with bound $B_\Psi$.
For those $g\in\h$, for which $\sum | \<g,\phi_n\>|^2 =\infty$ or $g=0$, clearly the lower frame condition holds. 
Now let $g\in \h$ be such that $\sum | \<g,\phi_n\>|^2<\infty$ and $g\neq 0$. For every $f\in\h$,
 $$
\mid \<Mf,g\>\mid
\leq 
\left(\sum | \<f,\psi_n\> |^2 \right)^{1/2}
\left(\sum | \<\phi_n,g\> |^2 \right)^{1/2}\leq 
\sqrt{B_\Psi} \|f\|  \left(\sum | \<\phi_n,g\> |^2 \right)^{1/2}.
$$
For $f=M^{-1}g$, it follows that
$$
\|g\|^2 \leq  
\sqrt{B_\Psi} \,  \, \|M^{-1}\| \,\|g\| \left(\sum | \<\phi_n,g\> |^2 \right)^{1/2}.
$$
Therefore, $\Phi$ satisfies the lower frame condition with bound $\frac{1}{B_\Psi\, \|M^{-1}\|^2}$. 

The case, when $\Phi$ is a Bessel sequence, can be shown in a similar way, using the inequality $|\<Mg,f\>| \leq \sqrt{B_\Phi} \|f\| (\sum | \<g,\psi_n\>|^2)^{1/2}$ applied with $f=Mg$. 

\v Second step: general $m$.

Apply the first step to the multiplier $M_{(1),m\Phi,\Psi}$ (resp. $M_{(1),\Phi,\overline{m} \,\Psi}$).

(ii) and (iii) follow easily from (i).

(iv) Let $\Psi$ and $\Phi$ be Bessel for $\h$ and $m\newin\ell^\infty$. Then $m\Phi$ and $m\Psi$ are also Bessel for $\h$.
The rest follows from (ii).
\ep

\v Note that Theorem \ref{c1}(i) generalizes one direction of \cite[Prop. 3.4]{casoleli1}, which states that if $\Phi$ is Bessel for $\h$ and $T_\Phi U_\Psi=I_\h$, then $\Psi$ fulfills the lower frame condition for $\h$. 
Furthermore, \cite[Lemma 5.6.2]{ole1}, which states that if two Bessel sequences fulfill $T_\Phi U_\Psi=I_\h$ then they are frames, is a special case of Theorem \ref{c1}(iv).  
To emphasize this, 
Theorem \ref{c1} can be rephrased into: 
\begin{cor} \
\begin{enumerate} 
\item Let $\Psi$ (resp. $\Phi$) be a Bessel sequence for $\h$ 
and let $T_\Phi U_\Psi$ be invertible on $\h$.
Then $\Phi$ (resp. $\Psi$) fulfills the lower frame condition for $\h$.

\item Let $\Psi$ and $\Phi$ be Bessel sequences for $\h$ such that $T_\Phi U_\Psi$ is invertible on $\h$. 
Then $\Psi$ and $\Phi$ are frames for $\h$.
\end{enumerate}
\end{cor}

By Theorem \ref{c1}(iv), a Bessel multiplier $M_{m,\Phi,\Psi}$ with $m\in\ell^\infty$ can 
be invertible only if the Bessel sequences $\Phi$ and $\Psi$ are frames. Note that the boundedness of $m$ is essential for this statement - if $m\notin\ell^\infty$, then a Bessel multiplier $M_{m,\Phi,\Psi}$ can be invertible also in cases when the Bessel sequences $\Phi$ and $\Psi$ are not frames. Examples are given in \cite{BStable09}.

\subsubsection{One of the sequences $\Phi$ and $\Psi$ is Bessel, which is not a frame}

Let $\Phi$ be $NBB$ Bessel for $\h$, which is not a frame for $\h$. 
In this case $M_{m,\Phi,\Psi}$ (resp. $M_{m,\Psi,\Phi}$) can not be both unconditionally convergent on $\h$ and invertible on $\h$ (see Propositions \ref{lem32}(i), \ref{uncnew} and Theorem \ref{c1}).

Let $\Phi$ be non-$NBB$ Bessel for $\h$, which is not a frame for $\h$. By Theorem \ref{c1}, if $m\Psi$ is Bessel for $\h$, then the multiplier $M_{m,\Phi,\Psi}$ (resp. $M_{m,\Psi,\Phi}$) can not be invertible on $\h$. If $m\Psi$ is non-Bessel for $\h$, then both cases (invertible multiplier and non-invertible multiplier) are possible - for example, $M_{(1),(\frac{1}{n}e_n),(ne_n)}=I_\h$ is invertible on $\h$ with unconditional convergence on $\h$ and $M_{(1),(\frac{1}{n^2}e_n),(ne_n)}$ is unconditionally convergent but not invertible on $\h$, see Example \ref{exnonsurj2}.

\subsubsection{One of the sequences $\Phi$ and $\Psi$ is a frame}

Let $\Phi$ be a $NBB$ frame for $\h$. In this case unconditional convergence and invertibility of $M_{m,\Phi,\Psi}$ (resp. $M_{m,\Psi,\Phi}$) on $\h$ require $m\Psi$ to be a frame for $\h$ (see Propositions \ref{lem32}(i), \ref{uncnew} and Theorem \ref{c1}(ii)). 
In the other direction, 
when $m\Psi$ is a frame for $\h$, both invertibility and non-invertibility of multipliers is possible - examples with overcomplete frames are given in Example \ref{exof}; the Riesz basis case is completely characterized in Section \ref{rb}.

Let $\Phi$ be a non-$NBB$ frame for $\h$ (hence, not a Riesz basis). 
If $m\Psi$ is a Riesz basis for $\h$ (resp. Bessel for $\h$ which is not a frame for $\h$), then $M_{m,\Phi,\Psi}$ can not be invertible on $\h$, see   Corollary \ref{rc3}(ii) (resp. Theorem \ref{c1}(ii)). If $m\Psi$ is non-Bessel or an overcomplete frame for $\h$, then both invertibility and non-invertibility of multipliers is possible. Many examples with consideration of $m$ - $SN$, $m\in\ell^\infty$, $m\notin\ell^\infty$, can be found in \cite{BStable09}.

We continue with sufficient conditions for invertibility of  $M_{m,\Phi,\Psi}$  and $M_{m,\Psi,\Phi}$.

\begin{prop} \label{gphi}
Let $\Phi$ be a frame for $\h$, $G:\h\to\h$ be a bounded bijective operator and $\psi_n=G\phi_n$, $\forall n$, i.e. $\Phi$ and $\Psi$ are equivalent frames
\footnote{For a treatise of equivalence of frames see \cite{antoin2} for the continuous and \cite{Casaz1} for the discrete setting.}. 
Let $m$ be positive (resp. negative) and semi-normalized.
Then $\Psi$ is a frame for $\h$, the multipliers 
$M_{m,\Phi,\Psi}$ and $M_{m,\Psi,\Phi}$ are invertible on $\h$ and
\begin{equation} \label{mginv} M_{m,\Phi,\Psi}^{-1}=
\left\{\begin{array}{ll}
(G^{-1})^* S^{-1}_{(\sqrt{m_n}\,\phi_n)}, \ &\mbox{when} \ m_n>0,\forall n,\\
-(G^{-1})^* S^{-1}_{(\sqrt{|m_n|}\,\phi_n)}, \ &\mbox{when} \ m_n<0,\forall n.
\end{array}
\right.
\end{equation}
\begin{equation} \label{mginv2} M_{m,\Psi,\Phi}^{-1}=
\left\{\begin{array}{ll}
 S^{-1}_{(\sqrt{m_n}\,\phi_n)} G^{-1}, \ &\mbox{when} \ m_n>0,\forall n,\\
- S^{-1}_{(\sqrt{|m_n|}\,\phi_n)} G^{-1}, \ &\mbox{when} \ m_n<0,\forall n.
\end{array}
\right.
\end{equation}
\end{prop}
\bp
By \cite[Corollary 5.3.2]{ole1}, $\Psi$ is a frame for $\h$.
For every $f\in \h$ we have $M_{m,\Phi,\Psi}f=
\sum m_n\<G^*f,\phi_n\>\phi_n =M_{m,\Phi,\Phi}G^*f.$
By \cite[Lemma 4.4]{xxljpa1}, $M_{m,\Phi,\Phi}$ is invertible
and
$$M_{m,\Phi,\Phi}^{-1}=
\left\{\begin{array}{ll}
 S^{-1}_{(\sqrt{m_n}\,\phi_n)}, \ &\mbox{when} \ m_n>0, \forall n,\\
- S^{-1}_{(\sqrt{\mid m_n\mid }\,\phi_n)}, \ &\mbox{when} \ m_n<0, \forall n.
\end{array}
\right.$$
Therefore,
 $M_{m,\Phi,\Psi}=M_{m,\Phi,\Phi}G^*$ is invertible and Equation (\ref{mginv}) holds. By similar arguments, Equation (\ref{mginv2}) holds.
\ep

\v 
As stated in the introduction and seen in some statements above, the invertibility of multipliers is related to the topic of dual frames.
Note that the above proposition covers the case when $\Psi$ is the canonical dual of $\Phi$ and does not cover any other dual frame of $\Phi$. Indeed, if $\Psi=(G\phi_n)$ is a dual frame of $\Phi$ for some bounded operator $G$, then $G$ must coincide with $S_\Phi^{-1}$, and thus, $\Psi$ must be the canonical dual of $\Phi$, see \cite[pp.19-20]{hala00}. For other duals, the following statement can be used.

\begin{prop} \label{p1}
Let $\Phi$ be a frame for $\h$ 
and let $\Phi^d=(\phi^d_n)$ be a dual frame of $\Phi$. 
Let $0\leq \lambda <\frac{1}{\sqrt{B_\Phi B_{\Phi^d}}}$
and let $(m_n)$ be such that $1-\lambda \leq m_n\leq 1+\lambda$, $\forall n\in\mn$. 
Then $M_{m,\Phi,\Phi^d}$ and $M_{m,\Phi^d,\Phi}$ are invertible on $\h$ and  \begin{equation}\label{edno2} \frac{1}{1+\lambda \sqrt{B_\Phi B_{\Phi^d}}}\|h\|\leq \|M^{-1} h\|\leq \frac{1}{1-\lambda \sqrt{B_\Phi B_{\Phi^d}}}\|h\|, \ \forall h\in\h,\end{equation}
where $M$ denotes any one of $M_{m,\Phi,\Phi^d}$ and $M_{m,\Phi^d,\Phi}$.
Moreover, $$M_{m,\Phi,\Phi^d}^{-1}=
\sum_{k=0}^\infty (M_{(1-m_n),\Phi,\Phi^d})^k \ \mbox{and} \ M_{m,\Phi^d,\Phi}^{-1}=\sum_{k=0}^\infty (M_{(1-m_n),\Phi^d,\Phi})^k.$$

\end{prop}
\bp 
The case $\lambda=0$ is trivial - in this case 
$M_{m,\Phi,\Phi^d}=M_{m,\Phi^d,\Phi}=M_{(1),\Phi,\Phi^d}=I_\h$. Now consider the case $\lambda > 0$.
First note that $B_{\Phi} B_{\Phi^d} \geq 1$. Indeed, for every $h\in\h$, $h\neq 0$, one has $\|h\|=\|T_{\Phi^d} U_\Phi h\|\leq \|T_{\Phi^d}\|\, \| U_\Phi\|\,\| h\|\leq \sqrt{B_{\Phi^d} B_\Phi}\,\|h\|$. Hence, $\lambda<1$ and $m$ is positive.
By assumption, $(m_n-1)\in\ell^\infty$. 
For every $f\in \h$, 
\begin{eqnarray*}
 \|M_{m,\Phi,\Phi^d}f-f\|&=&
\| \sum m_n \< f, \phi^d_n \> \phi_n - \sum \< f, \phi^d_n\>\phi_n \|=
\|M_{(m_n-1), \Phi, \Phi^d}f\|\\
&\leq& \lambda \sqrt{B_\Phi B_{\Phi^d}} \|f\|.
\end{eqnarray*}
  Since $\lambda\sqrt{B_\Phi B_{\Phi^d}}<1$, Proposition \ref{ginvnew} with $F=I_\h$ implies that the 
   multiplier $M_{m,\Phi, \Phi^d}$ is invertible on $\h$, (\ref{edno2}) holds and 
$$M_{m,\Phi,\Phi^d}^{-1}=\sum_{k=0}^\infty (I_\h-M_{m,\Phi,\Phi^d})^k=\sum_{k=0}^\infty (M_{(1-m_n),\Phi,\Phi^d})^k.$$
Since $\Phi$ is a dual of $\Phi^d$, the conclusions for $M_{m,\Phi^d,\Phi}$ follow directly from what is already proved. \ep

Let $\Phi$ be a frame for $\h$ and let $m$ be positive and semi-normalized. By Proposition \ref{gphi}, the multiplier $M_{m,\Phi,\widetilde{\Phi}}$ is invertible on $\h$ and $M_{m,\Phi,\widetilde{\Phi}}^{-1}=(S_\Phi)^* S^{-1}_{(\sqrt{m_n}\phi_n)}$. As a consequence of Proposition \ref{p1}, for certain weights $m$, the inverse $M_{m,\Phi,\widetilde{\Phi}}^{-1}$ can be expressed more simply via the use of $M_{m,\Phi,\widetilde{\Phi}}$: 

\begin{cor}\label{cp1}
Let $\Phi$ be a frame for $\h$ and $\widetilde{\Phi}$ be the canonical dual of $\Phi$. Let $0\leq \lambda <\sqrt\frac{A_\Phi}{B_\Phi}$ and let $\seq[m]$ be such that $1-\lambda \leq m_n\leq 1+\lambda$, $\forall n\in\mn$.
Then $M_{m,\Phi,\widetilde{\Phi}}$ and $M_{m,\widetilde{\Phi},\Phi}$ are invertible on $\h$ and
\begin{equation*} \frac{1}{1+\lambda \sqrt{B_\Phi/A_\Phi}}\|h\|\leq \|M^{-1}h\|\leq \frac{1}{1-\lambda \sqrt{B_\Phi/A_\Phi}}\|h\|, \ \forall h\in\h,\end{equation*}
where $M$ denotes any one of $M_{m,\Phi,\widetilde{\Phi}}$ and $M_{m,\widetilde{\Phi},\Phi}$.
Moreover,
$$M_{m,\Phi,\widetilde{\Phi}}^{-1}=
\sum_{k=0}^\infty (M_{(1-m_n),\Phi,\widetilde{\Phi}})^k \ \mbox{and} \ M_{m,\widetilde{\Phi},\Phi}^{-1}=\sum_{k=0}^\infty (M_{(1-m_n),\widetilde{\Phi},\Phi})^k.$$
\end{cor}
\bp
The assertion follows from Proposition \ref{p1}, because $\frac{1}{A_\Phi}$ is an upper frame bound for $\widetilde{\Phi}$.
\ep

Note that the bound for $\lambda$ in Corollary \ref{cp1} (resp. Proposition \ref{p1}) is sharp - if the assumptions hold with $\lambda =\sqrt\frac{A_\Phi}{B_\Phi}$ (resp. $\lambda=\sqrt{B_\Phi B_{\Phi^d}} $), then the multiplier might be non-invertible on $\h$ (see Example \ref{p1ex}). 

\v For the next assertion, we will use the property that, when $\Psi$ and $\Phi$ are Bessel sequences, then $\Psi \pm \Phi$ is a Bessel sequence with optimal bound $B_{\Psi \pm \Phi}^{opt}\leq (\sqrt{B_\Psi}+\sqrt{B_\Phi})^2$. Recall that if $\Phi$ is a frame for $\h$ and $m$ is positive (resp. negative) and semi-normalized, then $M_{m,\Phi,\Phi}$ is invertible on $\h$ (see \cite[Lemma 4.4]{xxljpa1}). Below, we generalize this statement allowing different sequences $\Phi$ and $\Psi$ in the multiplier. 

\begin{prop} \label{mpos}  
 Let $\Phi$ be a frame for $\h$. 
  Assume that $\Psi - \Phi$ is a Bessel sequence for $\h$ with bound $B_{\Psi-\Phi}< \frac{A_\Phi^2}{B_\Phi}$.
 For every positive (or negative) semi-normalized sequence $m$, satisfying
\begin{equation}\label{usl}
0< a \leq  \mid m_n \mid \leq b,\, \forall n, \  \mbox{and} \ \ \frac{b}{a}< \frac{A_\Phi}{\sqrt{B_{\Psi - \Phi} B_{\Phi}}},
\end{equation}
it follows that $\Psi$ is a frame for $\h$, 
the multipliers $M_{m,\Phi,\Psi}$ and $M_{m,\Psi,\Phi}$ are invertible on $\h$ and 
\begin{equation} \label{mnorm} \frac{1}{b B_\Phi + b\sqrt{B_\Phi B_{\Psi - \Phi}} } \|h\|\leq \|M^{-1} h\| \leq \frac{1}{aA_\Phi - b \sqrt{B_\Phi B_{\Psi - \Phi}} } \|h\|,\end{equation}
\begin{equation*}M^{-1}=
\left\{\begin{array}{ll}
\sum_{k=0}^\infty [S_{(\sqrt{m_n}\phi_n)}^{-1}(S_{(\sqrt{m_n}\phi_n)}-M)]^k S_{(\sqrt{m_n}\phi_n)}^{-1}, \ &\mbox{if} \ m_n>0,\forall n,\\
-\sum_{k=0}^\infty [S_{(\sqrt{|m_n|}\phi_n)}^{-1}(S_{(\sqrt{|m_n|}\phi_n)}+M)]^k S_{(\sqrt{|m_n|}\phi_n)}^{-1},  \ &\mbox{if} \ m_n<0, \forall n,
\end{array}
\right.\end{equation*}
where $M$ denotes any one of $M_{m,\Phi,\Psi}$ and $M_{m,\Psi,\Phi}$.
 \end{prop}
\bp 
Assume that $m$ is positive.
First note that we need $B_{\Psi-\Phi}< \frac{A_\Phi^2}{B_\Phi}$ in order to be able to fulfill (\ref{usl}). 
Assume now that (\ref{usl}) holds, and note that 
 $(\sqrt{m_n}\phi_n)$ is a frame for $\h$ with lower bound $a A_\Phi$ and upper bound $b B_\Phi$ (see \cite[Lemma 4.3]{xxljpa1}), and thus $\|S^{-1}_{(\sqrt{m_n}\phi_n)}\|\leq \frac{1}{aA_\Phi}$ (see Proposition \ref{Cprop}(iii)). 
 Since $\Phi$ is a frame for $\h$ and $B_{\Psi-\Phi}< \frac{A_\Phi^2}{B_\Phi}\leq A_\Phi$, it follows from \cite[Corollary 15.1.5]{ole1} that $\Psi$ is also a frame for $\h$. 
For every $h\in\h$,
$$
 \|M_{m,\Phi,\Psi}h-S_{(\sqrt{m_n}\phi_n)}h\|=\|M_{m,\Phi,\Psi-\Phi}\, h\| 
\leq b \cdot \sqrt{B_{\Psi-\Phi} B_{\Phi}}\, \|h\|.
$$
Because
$b \sqrt{B_{\Phi-\Psi} B_{\Phi}}< a A_{\Phi}\leq \frac{1}{\|S_{(\sqrt{m_n}\phi_n)}^{-1}\|} $, Proposition \ref{ginvnew} implies the invertibility of $M_{m,\Phi,\Psi}$, the representation for $M_{m,\Phi,\Psi}^{-1}$ and the upper inequality in (\ref{mnorm}). 
Moreover, for every $h\in\h$, 
\begin{eqnarray*}
\|M_{m,\Phi,\Psi}^{-1}h\| &\geq &  \frac{1}{\|M_{m,\Phi,\Psi}\|} \|h\|\geq \frac{1}{b\cdot \sqrt{B_\Phi}\sqrt{B^{opt}_{\Psi}}} \|h\| \\
&\geq& \frac{1}{b\cdot \sqrt{B_\Phi}(\sqrt{B_\Phi} + \sqrt{B_{\Psi-\Phi}})} \|h\|.
\end{eqnarray*}
We can derive a lower bound from Proposition \ref{ginvnew}, namely, $\frac{aA_\Phi}{bB_\Phi} \,  \frac{1}{a A_\Phi + b \cdot \sqrt{B_{\Phi-\Psi} B_{\Phi} } } $. However, the bound 
$\frac{1}{b(B_\Phi + \sqrt{B_\Phi B_{\Psi-\Phi}})}$ is sharper.

An analogous proof can be used for the invertibility of $M_{m,\Psi,\Phi}$ and the conclusions for $M_{m,\Psi,\Phi}^{-1}$.

If $m$ is negative, apply what is already proved to the multiplier $M_{-m,\Phi,\Psi}$.
\ep

\v Note that the bound for $b/a$ in (\ref{usl}) is sharp. If $\frac{\sup_n |m_n|}{\inf_n |m_n|}=\frac{A_\Phi}{\sqrt{B_{\Psi - \Phi} B_{\Phi}}}$, the multiplier $M_{m,\Phi,\Psi}$ can be non-invertible on $\h$ (see Example \ref{noninvex} with $m=(1)$ and 
$k\in[0,2]$).

\begin{prop} \label{p4}
Let $\Phi$ be a frame for $\h$. 
Assume that
\begin{itemize}
 \item[$\mathcal{P}_1${\rm :} ] 
$ \exists \,\mu\in[0,\frac{A_\Phi^2}{B_{\Phi}}) \ \mbox{such that} $
 \mbox{$\sum |\<f, m_n\psi_n-\phi_n \>|^2 \leq \mu \|f\|^2$, \, $\forall$ $f \in \h$.}
 \end{itemize}
Then $m\Psi$ is a frame for $\h$, the multipliers $M_{\overline{m},\Phi,\Psi}$ and $M_{m,\Psi,\Phi}$ are invertible on $\h$ and 
     \begin{equation}\label{multineq}
     \frac{1}{B_\Phi+\sqrt{\mu B_\Phi}}\|h\|
   \leq \|M^{-1}h\|
   \leq \frac{1}{A_\Phi-\sqrt{\mu B_\Phi}} \|h\|, \ \forall h\in\h,
   \end{equation}
\begin{equation}\label{minv}
M^{-1}=\sum_{k=0}^\infty [S_\Phi^{-1}(S_\Phi-M)]^k S_\Phi^{-1} \end{equation}
where $M$ denotes any one of $M_{\overline{m},\Phi,\Psi}$ and $M_{m,\Psi,\Phi}$.

As a consequence, if $m$ is semi-normalized, then $\Psi$ is also a frame for $\h$.

\end{prop}
\bp
If $\mathcal{P}_1$ holds with $\mu=0$, then $m\Psi=\Phi$ and thus, $M_{\overline{m},\Phi,\Psi}=M_{m,\Psi,\Phi}=S_\Phi$. 
Now apply Proposition \ref{Cprop}(iii).

Let $\mathcal{P}_1$ hold with $\mu\neq 0$. Apply Proposition \ref{mpos} to the multipliers $M_{(1),\Phi,m\Psi}$ and $M_{(1),m\Psi,\Phi}$.
\ep

\v Note that $\mathcal{P}_1$ is equivalent to the following two conditions:

\begin{itemize}
\item[$\bullet$] $ \exists \,\nu\in[0,\frac{A_\Phi}{\sqrt{B_{\Phi}}}) \ \mbox{such that} $

$\| \sum c_n (m_n\psi_n - \phi_n)\|\leq \nu \|\seq[c]\|_2 \ \mbox{for all finite scalar sequences}\ \seq[c]$ (and hence for all $(c_n)\in\ell^2$);
\item[$\bullet$]
$ \exists \,\mu\in[0,\frac{A_\Phi^2}{B_{\Phi}}) \ \mbox{such that} $

 \mbox{$\sum |\<f, m_n\psi_n-\phi_n \>|^2 \leq \mu \|f\|^2$, for all $f$ in a dense subset of $\h$.}
  \end{itemize}
Indeed, the case $\mathcal{P}_1$ \& ($\mu=0$) is trivial; for the case $\mathcal{P}_1$ \& ($\mu\neq 0$)  
see Proposition \ref{Cprop}(ii) and \cite[Lemma 3.2.6]{ole1}.

\v The bound for $\mu$ in Proposition \ref{p4} is sharp - if 
$\mu\geq \frac{A_\Phi^2}{B_{\Phi}}$, 
then the multiplier might be non-invertible (see Example \ref{noninvex}).
The multiplier might be invertible for any value of $\mu$ (see Example \ref{inv22}).

Although Proposition \ref{p4} is proved based on Proposition \ref{mpos}, the two propositions cover different classes of sequences $\Phi,\Psi,$ and $m$, even when $m$ is restricted to be real. 
 Example \ref{8yes7no} shows a case when Proposition \ref{p4} applies, but Proposition \ref{mpos} does not apply. For a case where Proposition \ref{mpos} can be used, but Proposition \ref{p4} can not be used, see Example \ref{exnew}.

By Proposition \ref{p4}, when $m\Psi$ is a perturbation of $\Phi$, the inverse operator $M_{m,\Psi,\Phi}^{-1}$ is given by (\ref{minv}). Simpler representation for $M_{m,\Psi,\Phi}^{-1}$ can be obtained if $m\Psi$ is a perturbation of a dual frame of $\Phi$:

\begin{prop} \label{p3}
 Let $\Phi$ be a frame for $\h$.
Assume that 
\begin{itemize}
 \item[$\mathcal{P}_2${\rm :} ] 
$ \exists \,\mu\in[0,\frac{1}{B_{\Phi}}) \ \mbox{such that} $
 \mbox{$\sum |\<f, m_n\psi_n-\phi^d_n \>|^2 \leq \mu \|f\|^2$, \ $\forall$ $f\in\h$,}
 \end{itemize}
\noindent for some dual frame $\Phi^d=(\phi^d_n)$ of $\Phi$.  Then $m\Psi$ is a frame for $\h$, the multipliers $M_{\overline{m},\Phi,\Psi}$ and $M_{m,\Psi,\Phi}$ are invertible on $\h$ and
\begin{equation*} \frac{1}{1+ \sqrt{\mu B_\Phi}}\|h\|\leq \|M^{-1}h\|\leq \frac{1}{1- \sqrt{\mu B_\Phi}}\|h\|, \ \forall h\in\h,\end{equation*}
 $$M^{-1}=\sum_{k=0}^\infty (I_\h-M)^k,$$
 where $M$ denotes any one of $M_{\overline{m},\Phi,\Psi}$ and $M_{m,\Psi,\Phi}$.
 
 As a consequence, if $m$ is semi-normalized, then $\Psi$ is also a frame for $\h$.
 \end{prop}
\bp 
Since $\Phi^d$ is a dual frame of $\Phi$, the number $\frac{1}{B_\Phi}$ is a lower bound for $\Phi^d$ (see the proof of \cite[Proposition 3.4]{casoleli1}). Since $\mu$ is smaller then the lower bound $\frac{1}{B_\Phi}$ of $\Phi^d$, 
it follows from \cite[Corollary 15.1.5]{ole1} that $m\Psi$ is a frame for $\h$.
Therefore, $M_{\overline{m},\Phi,\Psi}=M_{(1),\Phi,m\Psi}$ is well defined. 
For every $f\in\h$,
$$
\|M_{\overline{m},\Phi,\Psi}f-f\|=\|M_{(1),(\phi_n),(m_n\psi_n-\phi^d_n)} f\| 
\leq \sqrt{\mu B_{\Phi}}\|f\|
$$
and similarly, $ \|M_{m,\Psi,\Phi}f-f\|\leq\sqrt{\mu B_{\Phi}}\|f\|.$
Since $\mu B_{\Phi}\in [0,1)$, one can apply Proposition \ref{ginvnew} and this concludes the proof.
\ep

\v Similar to the case with $\mathcal{P}_1$, one can  list conditions equivalent to $\mathcal{P}_2$.

\v The bound for $\mu$ in Proposition \ref{p3} is sharp - if 
$\mu\geq \frac{1}{B_{\Phi}}$,
then 
the multiplier might be non-invertible (see Example \ref{noninvex}).
Note that the multiplier might be invertible for any value of $\mu$ (see Example \ref{inv22}).

\v Note that Propositions \ref{p4} and \ref{p3} do not cover the same classes of sequences. For a case when Proposition \ref{p4} applies (resp. does not apply) and Proposition \ref{p3} does not apply (resp. applies) see Example \ref{exdual} (resp. \ref{exdual2}).
Let $\Phi$ be a Parseval frame, i.e. a tight frame with $A=B=1$. In this case, the frame is self-dual, and both Propositions \ref{p4} and \ref{p3} can be applied.

Propositions \ref{mpos}  and \ref{p3} do not cover the same classes of sequences either. Example \ref{7yes9no} (resp. \ref{9yes7no}) shows a case when Proposition \ref{mpos} applies (resp. does not apply), but Proposition \ref{p3} does not apply (resp. applies).

\subsubsection{One of the sequences $\Phi$ and $\Psi$ is a Riesz basis} \label{rb}

For two Riesz bases and a semi-normalized symbol, the multipliers are always invertible \cite[Prop. 7.7]{xxlmult1}.  
If $\Phi$ is a Riesz basis for $\h$, $m$ is real and $SN$, and $\Psi$ is a frame for $\h$, then the multiplier $M_{m,\Phi,\Psi}$ (resp. $M_{m,\Psi,\Phi}$) is invertible on $\h$ if and only if $\Psi$ is a Riesz basis for $\h$  \cite[Prop. 4.2]{stoevxxl09}.
What can be said about the cases, when one of the sequences has the Riesz property, $m$ is complex and not necessarily semi-normalized, and $\Psi$ is not necessarily a frame? The answer is given in the following assertion. 

\begin{thm} \label{r123}
Let $\Phi$ be a Riesz basis for $\h$. 
\begin{itemize}
\item[{\rm (a)}]  Let $m$ be $SN$. 
\begin{itemize}
\item[{\rm (a1)}] 
If $\Psi$ is non-Bessel for $\h$, then both $M_{m,\Phi,\Psi}$ and $M_{m,\Psi,\Phi}$ are not well defined.

\item[{\rm (a2)}] If $\Psi$ is Bessel for $\h$, then
$M_{m,\Phi,\Psi}$ (resp. $M_{m,\Psi,\Phi}$) is invertible on $\h$ if and only if $\Psi$ is a Riesz basis for $\h$. In the cases of invertibility, 
$M_{m,\Phi,\Psi}^{-1}=M_{(\frac{1}{m_n}),\widetilde{\Psi},\widetilde{\Phi}}$.
\end{itemize}

\item[{\rm (b)}] Let $m$ be non-$NBB$ and $m\in\ell^\infty$. 

\begin{itemize}
\item[{\rm (b1)}] 
If $\Psi$ is $NBA$ non-Bessel for $\h$, then $M_{m,\Phi,\Psi}$ (resp. $M_{m,\Psi,\Phi}$) can either be well defined or not, but can never be invertible on $\h$.
\item[{\rm (b2)}] If $\Psi$ is non-$NBA$, non-$NBB$, and non-Bessel for $\h$, then $M_{m,\Phi,\Psi}$ (resp. $M_{m,\Psi,\Phi}$) can either be well defined or not, but can never be invertible on $\h$.

\item[{\rm (b3)}] 
If $\Psi$ is non-$NBA$, $NBB$, and non-Bessel for $\h$, then for $M_{m,\Phi,\Psi}$ and $M_{m,\Psi,\Phi}$, all the three feasible combinations of invertibility and well-definedness are possible: they can be invertible on $\h$, can be well defined and non-invertible on $\h$, or can be not well defined.

\item[{\rm (b4)}]  
If $\Psi$ is Bessel for $\h$, then $M_{m,\Phi,\Psi}$ (resp. $M_{m,\Psi,\Phi}$) is well defined on $\h$, but not invertible on $\h$.
\end{itemize}

\item[{\rm (c)}] Let $m$ be $NBB$ and $m\notin\ell^\infty$. 
\begin{itemize}
\item[{\rm (c1)}] 
If $\Psi$ is  non-Bessel for $\h$ or $NBB$, 
then both $M_{m,\Phi,\Psi}$ and $M_{m,\Psi,\Phi}$ are not well defined. 
\item[{\rm (c2)}] Let $\Psi$ be non-$NBB$ Bessel for $\h$, which is not a frame for $\h$. Then for $M_{m,\Phi,\Psi}$ and $M_{m,\Psi,\Phi}$ all feasible combinations of invertibility and well-definedness are possible.
\item[{\rm (c3)}] Let $\Psi$ be a non-$NBB$ frame for $\h$. Then $M_{m,\Phi,\Psi}$ (resp. $M_{m,\Psi,\Phi}$) can either be well defined or not, but can never be invertible on $\h$.

\end{itemize}

\item[{\rm (d)}] Let $m$ be non-$NBB$ and $m\notin\ell^\infty$. 
\begin{itemize}
\item[{\rm (d1)}] If $\Psi$ is $NBB$, then both $M_{m,\Phi,\Psi}$ and $M_{m,\Psi,\Phi}$ are not well defined. 
\item[{\rm (d2)}] 
If $\Psi$ is non-$NBB$ non-$NBA$ non-Bessel for $\h$, then for $M_{m,\Phi,\Psi}$ (resp. $M_{m,\Psi,\Phi}$) all feasible combinations of invertibility and well-definedness are possible.
\item[{\rm (d3)}] If $\Psi$ is non-$NBB$ and $NBA$, then $M_{m,\Phi,\Psi}$ (resp. $M_{m,\Psi,\Phi}$) can either be well defined or not, but can never be invertible on $\h$. 
\end{itemize}

\end{itemize}
\end{thm}
\bp
(a1) follows from Proposition \ref{lem32}(ii), because in the case when $m$ is $SN$, the sequence $m\Psi$ is Bessel for $\h$ if and only if $\Psi$ is Bessel for $\h$.

(a2) If $\Psi$ is Bessel for $\h$, which is not a frame for $\h$, then Theorem \ref{c1}(iv) implies that both  $M_{m,\Phi,\Psi}$ and $M_{m,\Psi,\Phi}$ are not invertible on $\h$. 

Let $\Psi$ be an overcomplete frame for $\h$. 
If $m$ is real, it is proved in \cite{stoevxxl09} that $M_{m,\Phi,\Psi}$ and $M_{m,\Psi,\Phi}$ are not invertible on $\h$. Almost the same proof can be used in the case when $m$ is complex, but for the sake of completeness we include a proof here.
First observe that the sequences $m\Psi$ and $\overline{m}\,\Psi$ are also overcomplete frames for $\h$, 
which implies that $\range(U_{\overline{m}\,\Psi})\subsetneqq \ell^2$ and $T_{m\Psi}$ is not injective.
Since $M_{m,\Phi,\Psi}=T_{\Phi} U_{\overline{m}\,\Psi}$ and $T_{\Phi}$ is a bijection of $\ell^2$ onto $\h$, it follows that $M_{m,\Phi,\Psi}$ is not surjective. 
Since $U_\Phi$ is a bijection of $\h$ onto $\ell^2$, 
it follows that $M_{m,\Psi,\Phi}=T_{m\Psi} U_{\Phi}$ is not injective.

When $\Psi$ is a Riesz basis for $\h$, the invertibility of $M_{m,\Phi,\Psi}$ and the representation for the inverse are proved in {\rm \cite[Prop. 7.7]{xxlmult1}}. 
 
 For the proof of (b)-(d), first recall that $M_{m,\Phi,\Psi}$ (resp. $M_{m,\Psi,\Phi}$) is well defined on $\h$ if and only if $m\Psi$ is Bessel for $\h$ (see Proposition \ref{lem32}(ii)).
 
 (b1) As an example of a well defined multiplier, consider the sequences $\Phi=(e_n)$, $\Psi=(e_1,e_2,e_1,e_3,e_1,e_4,\ldots)$, and $m=(\frac{1}{2},1,\frac{1}{2^2},1,\frac{1}{2^3},1,\ldots)$. Since $m\Psi$ is Bessel for $\h$, both $M_{m,\Phi,\Psi}$ and $M_{m,\Psi,\Phi}$ are well defined on $\h$. 
 Now consider the sequences $\Phi=(e_n)$, $\Psi=(e_1,e_2,e_1,e_3,e_1,e_4,\ldots)$, and $m=(1,\frac{1}{2},1,\frac{1}{3},1,\frac{1}{4},\ldots)$ - in this case both $M_{m,\Phi,\Psi}$ and $M_{m,\Psi,\Phi}$ are not well defined, because $m\Psi$ is not Bessel for $\h$.

Now assume that $M_{m,\Phi,\Psi}$ (resp. $M_{m,\Psi,\Phi}$) is well defined on $\h$  
and thus $m\Psi$ is Bessel for $\h$. 
By (a) and Proposition \ref{mmbar} (resp. (a)), the multiplier $M_{(1),\Phi,\overline{m}\,\Psi}$ (resp. $M_{(1),m\Psi,\Phi}$) is invertible on $\h$ if and only if $m\Psi$ is a Riesz basis for $\h$.
By Proposition \ref{ovfr}, the sequence  $m\Psi$ can not be a Riesz basis for $\h$ under the assumptions of (b1).

(b2) As an example of a well defined multiplier, consider the sequences  $\Phi=(e_n)$, $\Psi=(e_1,2e_2,\frac{1}{3}e_3,4e_4,\frac{1}{5}e_5,6e_6,\ldots)$, and $m=(1,\frac{1}{2},1,\frac{1}{4},1,\frac{1}{6},\ldots)$. Since $m\Psi$ is Bessel for $\h$, both $M_{m,\Phi,\Psi}$ and $M_{m,\Psi,\Phi}$ are well defined on $\h$. 
Now consider the same $\Phi$ and $\Psi$, and the sequence $\nu=(1,1,\frac{1}{3},1,\frac{1}{5},1,\ldots)$. Both $M_{\nu,\Phi,\Psi}$ and $M_{\nu,\Psi,\Phi}$ are not well defined, because $\nu\Psi$ is not Bessel for $\h$ (use Lemma \ref{lem31}).  

Now assume that $M_{m,\Phi,\Psi}$ (resp. $M_{m,\Psi,\Phi}$) is invertible on $\h$. It follows from (a) that $m\Psi$ is a Riesz basis for $\h$. Hence, there exists $a>0$ so that $\|m_n\psi_n\|\geq a$, $\forall n$. Since $m\in\ell^\infty$, the last inequality implies that $\Psi$ is $NBB$, which contradicts to the assumptions. 
 
(b3) As an example of invertible multipliers on $\h$, consider 
 $M_{(1/n),(e_n),(ne_n)}=M_{(1/n),(ne_n),(e_n)}=I_\h$. As an example for well-defined non-invertible multipliers, take  $M_{(1/n^2),(e_n),(ne_n)}=M_{(1/n^2),(ne_n),(e_n)}$ (see Example \ref{exnonsurj2}). 
For a case with multipliers which are not well defined, consider $M_{(\frac{1}{n}), (e_n),(n^2e_n)}$ and $M_{(\frac{1}{n}), (n^2e_n),(e_n)}$.

(b4) 
Proposition \ref{bmh} implies that $M_{m,\Phi,\Psi}$ (resp. $M_{m,\Psi,\Phi}$) is well defined on $\h$.
The non-invertibility of $M_{m,\Phi,\Psi}$ (resp. $M_{m,\Psi,\Phi}$) can be shown in an analogue way as in (b1).
 
 (c1) 
 By Proposition \ref{lem32}(ii) and Corollary \ref{corlem32}, well-definedness of $M_{m,\Phi,\Psi}$ (resp. $M_{m,\Psi,\Phi}$) requires $\Psi$ to be Bessel for $\h$.
 If $\Psi$ is $NBB$, the conclusion follows from Proposition \ref{lem32}(iii).
 
(c2) For a case with invertible multipliers look at $M_{(n),(e_n),(\frac{1}{n}e_n)}=M_{(n),(\frac{1}{n}e_n),(e_n)}=I_\h$. As an example of well defined non-invertible multipliers, take $M_{(n),(e_n),(\frac{1}{n^2}e_n)}=M_{(n),(\frac{1}{n^2}e_n),(e_n)}$, see Example \ref{exnonsurj2}. The multipliers $M_{(n^2),(e_n),(\frac{1}{n}e_n)}=M_{(n^2),(\frac{1}{n}e_n),(e_n)}$ are not well defined. 
  
(c3) Consider $\Phi=(e_n)$ and the sequence  $\Psi=(\frac{1}{2}e_1, e_2, \frac{1}{2^2}e_1, e_3, \frac{1}{2^3}e_1, e_4, \ldots)$, which  is non-$NBB$ frame for $\h$. For $m=(\sqrt{2},1,\sqrt{2^2},1,\sqrt{2^3},1,\ldots)$, the sequence $m\Psi$ is  Bessel for $\h$, which implies that both $M_{m,\Phi,\Psi}$ and $M_{m,\Psi,\Phi}$ are well defined on $\h$.
For $m=(2,1,2^2,1,2^3,1,\ldots)$, the sequence $m\Psi$ is not Bessel for $\h$, which implies that both $M_{m,\Phi,\Psi}$ and $M_{m,\Psi,\Phi}$ are not well defined. If $M_{m,\Phi,\Psi}$ (resp. $M_{m,\Psi,\Phi}$) is well defined on $\h$, the non-invertibility of $M_{m,\Phi,\Psi}$ (resp. $M_{m,\Psi,\Phi}$) can be shown in an analogue way as in (b1).

(d1) follows from Proposition \ref{lem32}(iii).

(d2) 
 Consider $\Phi=(e_n)$ and the sequence $\Psi=(e_1,\frac{1}{2^2}e_2, 3e_3, \frac{1}{4^2}e_4, 5e_5, \ldots)$. For $m=(1,2^2,\frac{1}{3},4^2,\frac{1}{5},\ldots)$, we have $M_{m,\Phi,\Psi}=M_{m,\Psi,\Phi}=I_\h$.
For $m=(1,2,\frac{1}{3^2},4,\frac{1}{5^2},\ldots)$, both multipliers $M_{m,\Phi,\Psi}$ and $M_{m,\Psi,\Phi}$ coincide with the non-invertible operator used in Example \ref{exnonsurj2}.
For $m=(1,2^3,\frac{1}{3},4^3,\frac{1}{5},\ldots)$, the sequence $m\Psi$ is not Bessel for $\h$, which implies that both $M_{m,\Phi,\Psi}$ and $M_{m,\Psi,\Phi}$ are not well defined.

(d3) Examples for the case \lq\lq $\Psi$ - non-NBB Bessel\rq\rq:

Let $\Phi=(e_n)$ and $\Psi=(e_1,\frac{1}{2}e_2,\frac{1}{3}e_3,\ldots)$. For $m=(1,2^2,\frac{1}{3},4^2,\frac{1}{5},\ldots)$, $m\Psi$ is not Bessel for $\h$, and thus both $M_{m,\Phi,\Psi}$ and $M_{m,\Psi,\Phi}$ are not well defined. 
 For $m=(1,2,\frac{1}{3},4,\frac{1}{5},\ldots)$, $m\Psi$ is Bessel for $\h$,  and thus both $M_{m,\Phi,\Psi}$ and $M_{m,\Psi,\Phi}$ are well defined. 
 
 Examples for the case \lq\lq $\Psi$ - $NBA$ non-NBB non-Bessel\rq\rq:
 
 Let $\Phi=(e_n)$ and $\Psi=(e_1,\frac{1}{2}e_2,e_1,\frac{1}{3}e_3,e_1,\frac{1}{4}e_4\ldots)$. For $m=(\frac{1}{2}, 2, \frac{1}{2^2},3,\frac{1}{2^3},4,\ldots)$, the sequence $m\Psi$ is Bessel for $\h$ and thus, both $M_{m,\Phi,\Psi}$ and $M_{m,\Psi,\Phi}$ are well defined. 
 For $m=(\frac{1}{2}, 2^2, \frac{1}{2^2},3^2,\frac{1}{2^3},4^2,\ldots)$, $m\Psi$ is non-Bessel for $\h$  and thus, both $M_{m,\Phi,\Psi}$ and $M_{m,\Psi,\Phi}$ are not well defined. 
 
If $M_{m,\Phi,\Psi}$ (resp. $M_{m,\Psi,\Phi}$) is well defined on $\h$, the non-invertibility of $M_{m,\Phi,\Psi}$ (resp. $M_{m,\Psi,\Phi}$) can be shown in an analogue way as in (b1).
   \ep

As a consequence of the above detail assertion and having in mind Proposition \ref{mmbar}, we can summarize the possibilities for invertibility as follows:

\begin{cor} \label{rc3}
Let $\Phi$ be a Riesz basis for $\h$. Then $M_{m,\Phi,\Psi}$ (resp. $M_{m,\Psi,\Phi}$) is invertible on $\h$ if and only if $m\Psi$ is a Riesz basis for $\h$. 

Further, the following holds.
\begin{itemize}
\item[{\rm (i)}] If $\Psi$ is a Riesz basis for $\h$, then $M_{m,\Phi,\Psi}$ (resp. $M_{m,\Psi,\Phi}$) is invertible on $\h$ if and only if  $m$ is $SN$.
\item[{\rm (ii)}] If $m$ is $SN$, then $M_{m,\Phi,\Psi}$ (resp. $M_{m,\Psi,\Phi}$) is invertible on $\h$ if and only if $\Psi$ is a Riesz basis for $\h$. 
\item[{\rm (iii)}] If $m$ is not $SN$, then $M_{m,\Phi,\Psi}$ (resp. $M_{m,\Psi,\Phi}$) can be invertible on $\h$ only in the following cases:

$\bullet$ $\Psi$ is non-$NBB$ and Bessel for $\h$, which is not a frame for $\h$, and $m$ is $NBB$, but not in $\ell^\infty$;

$\bullet$ $\Psi$ is non-$NBA$, $NBB$, and non-Bessel for $\h$, $m$ is non-$NBB$ and $m\in\ell^\infty$;

$\bullet$ $\Psi$ is non-$NBA$, non-$NBB$, and non-Bessel for $\h$, $m$ is non-$NBB$ and $m\notin\ell^\infty$.

In the cases of invertibility, $M_{m,\Phi,\Psi}^{-1}=M_{(1),\widetilde{\overline{m}\,\Psi},\widetilde{\Phi}}$ and $M_{m,\Psi,\Phi}^{-1}=M_{(1),\widetilde{\Phi},\widetilde{m\,\Psi}}$. 
\end{itemize}
\end{cor}

\section{Examples}

In this section we list some examples, which we refer to throughout the paper.

\begin{ex}\label{exnonsurj2}
The operator $M:\h\to\h$ given by $M f=\sum \frac{1}{n}\<f,e_n\>e_n$ is injective, but not surjective - for example, the element $\sum \frac{1}{n}e_n\in\h$ does not belong to the range of $M$.
  \end{ex}

\begin{ex}\label{exof}
Invertible multiplier of two overcomplete frames:  

$\Phi=(e_1,e_1,e_2,e_3,e_4,\ldots)$, $\Psi=(\frac{1}{2}e_1,\frac{1}{2}e_1,e_2,e_3,e_4,\ldots)$,  $M_{(1),\Phi,\Psi}=I_\h$ with  unconditional convergence on $\h$.

Non-invertible multiplier of two overcomplete frames: 

$\Phi=(e_1,e_1,e_2,e_2,e_3,e_3,\ldots)$,  $\Psi=(e_1,e_1,e_2,e_3,e_4,\ldots)$, $M_{(1),\Phi,\Psi}$ is unconditionally convergent on $\h$, but not injective. 
\end{ex}

\subsection{Examples for the sharpness of the bounds:}
The bound for $\lambda$ in Proposition \ref{p1} is sharp:

\begin{ex} \label{p1ex} Let $\Phi=(e_n)$ and $m=(\frac{1}{n})$. The smallest possibility for $\lambda$ satisfying  $\sup_n \mid m_n - 1\mid \leq \lambda$ is $1$ and $1=1/\sqrt{B_\Phi B_{\widetilde{\Phi}}}=\sqrt{\frac{A_\Phi}{B_\Phi}}$. The multiplier $M_{m,\Phi,\widetilde{\Phi}}$ is injective, but not surjective (see Example \ref{exnonsurj2}).
\end{ex}

The bound for $\mu$ in Proposition \ref{p4} (resp. \ref{p3}) is sharp. If $\mathcal{P}_1$ (resp. $\mathcal{P}_2$) holds with $\mu\geq A_\Phi^2/B_\Phi$ (resp. $\mu\geq 1/B_\Phi$), then the multipliers $M_{\overline{m},\Phi,\Psi}$ and $M_{m,\Psi,\Phi}$ might be non-invertible on $\h$:

\begin{ex} \label{noninvex}
Let $\Phi=( e_n)$ and $m\Psi=(ke_1, \frac{1}{2}\,e_2,\frac{1}{3}\,e_3,\frac{1}{4}\,e_4, \ldots)$ for some number $k$. The unique dual frame of $\Phi$ is $\Phi^d=(e_n)$.
The multipliers $M_{(1),\Phi,m\Psi}$ and $M_{(1),m\Psi,\Phi}$ are non-invertible on $\h$ by Theorem \ref{c1}(iv), because  $\Phi$ is a frame for $\h$ and $m\Psi$ is Bessel for $\h$, which is not a frame for $\h$.
 
The sequence $m\Psi - \Phi=m\Psi-\Phi^d$ satisfies:
\begin{eqnarray*}
\sum |\<h,m_n\psi_n-\phi_n\>|^2&=&|k-1|^2 \, |\<h,e_1\>|^2 + \sum_{n=2}^\infty \left(\frac{1}{n}-1\right)^2 |\<h,e_n\>|^2.
\end{eqnarray*}
If $|k-1|\leq 1$, then 
$\sum |\<h,m_n\psi_n-\phi_n\>|^2\leq \|h\|^2,\ \forall h\in\h, \ $ and $$ \sum |\<e_i,m_n\psi_n-\phi_n\>|^2=\left(\frac{1}{i}-1\right)^2 \|e_i\|^2,\ \forall i=2,3,4,\ldots,$$
which implies that $m\Psi - \Phi$ is Bessel for $\h$ with optimal bound equal to $1$.

If $|k-1|> 1$, then 
$\sum |\<h,m_n\psi_n-\phi_n\>|^2\leq |k-1|^2 \, \|h\|^2, \ \forall h\in\h$ and
$$\sum |\<e_1,m_n\psi_n-\phi_n\>|^2=|k-1|^2 \, \|e_1\|^2,$$
which implies that $m\Psi - \Phi$ is Bessel for $\h$ with an optimal bound equal to $|k-1|^2$. 
Therefore, 
\begin{equation*} B^{opt}_{m\Psi-\Phi}=
\left\{\begin{array}{ll}
|k-1|^2>1=1/B_\Phi=A_\Phi^2/B_\Phi, \ &\mbox{when} \ |k-1|>1,\\
1=1/B_\Phi=A_\Phi^2/B_\Phi, \ &\mbox{when} \ |k-1|\leq 1,
\end{array}
\right.\end{equation*}
which shows that the example fulfills $\mathcal{P}_1$ (resp. $\mathcal{P}_2$) with any $\mu\geq A_\Phi^2/B_\Phi $ (resp. $\mu\geq 1/B_\Phi$). 
\end{ex}

\v Note that the multipliers $M_{\overline{m},\Phi,\Psi}$ and $M_{m,\Psi,\Phi}$ can be invertible with any value of $\mu$ in $\mathcal{P}_1$ (resp. $\mathcal{P}_2$) :

\begin{ex} \label{inv22} Let $\Phi=(e_n)$. The unique dual frame of $\Phi$ is $\Phi^d=(e_n)$.
\begin{itemize}
\item[{\rm (a)}]  Let $m\Psi=(ke_1, e_2, e_3, e_4, \ldots)$, where $k\neq 0$. 
The sequence $m\Psi-\Phi=m\Psi - \Phi^d$  
is Bessel for $\h$ with optimal bound $\mu=|k-1|^2\neq 1=A_\Phi^2/B_\Phi=1/B_\Phi$. 
The multipliers $M_{(1),\Phi,m\Psi}$ and $M_{(1),m\Psi,\Phi}$ are invertible on $\h$. 
\item[{\rm (b)}] Let $m\Psi=(2e_n)$. The sequence $m\Psi-\Phi=m\Psi - \Phi^d=(e_n)$  
is Bessel for $\h$ with optimal bound $\mu=1=A_\Phi^2/B_\Phi=1/B_\Phi$. 
The multipliers $M_{(1),\Phi,m\Psi}$ and $M_{(1),m\Psi,\Phi}$ are equal to $2 I_\h$ and thus, they are invertible on $\h$.
\end{itemize}
\end{ex}

\subsection{Independence of conditions for invertibility of multipliers}

{\bf Proposition \ref{p4} applies, Proposition \ref{p3} does not apply:}

\begin{ex} \label{exdual} Let $\Phi=( e_1, e_1, e_2, e_3, e_4, \ldots )$. Clearly, $\Phi$ is a frame for $\h$ with $A_\Phi^{opt}=1$, $B_\Phi^{opt}=2$. Take $m_n\psi_n=(k+1)\phi_n$, $\forall n$, where $k\in (0,\frac{1}{2})$. The sequence $m\Psi-\Phi=(k\phi_n)$ is Bessel for $\h$ with an optimal bound $B_{m\Psi-\Phi}^{opt}=k^2B^{opt}_\Phi=2k^2<\frac{1}{2}=\frac{(A^{opt}_\Phi)^2}{B^{opt}_\Phi}$. Thus, Proposition \ref{p4} implies the invertibility of  $M_{\overline{m},\Phi,\Psi}$ and $M_{m,\Psi,\Phi}$.

Now observe that the sequences $(h,e_1-h,e_2,e_3,e_4,\ldots)$, $h\in\h$, are precisely all the dual frames of $\Phi$.
Let $\Phi^d=(\phi^d_1, e_1-\phi^d_1,e_2,e_3,e_4,\ldots)$ be an arbitrary chosen dual frame of $\Phi$ and denote $\<e_1,\phi^d_1\>=x+iy$. 
 Consider the sequence
$$(m_n\psi_n - \phi^d_n)=((k+1)e_1 - \phi^d_1, k e_1 +\phi^d_1, ke_2, ke_3, ke_4,\ldots)$$
and observe that
\begin{eqnarray*}
\sum | \< e_1, m_n\psi_n - \phi^d_n \> |^2&=&| \<e_1, (k+1)e_1 - \phi_1^d\> |^2 + | \<e_1, ke_1 + \phi_1^d\> |^2 \\
&=&| k+1 - x-iy |^2 + | k + x + iy |^2 \\
&=& \left(2k^2 + 2x^2 + 2k -2x + 2y^2 + 1\right)\|e_1\|^2.
\end{eqnarray*}
 Assume that $m\Psi-\Phi^d$ is Bessel for $\h$ with bound $B_{m\Psi-\Phi^d}\leq \frac{1}{2}(=\frac{1}{B^{opt}_\Phi})$. Applying the Bessel inequality to the element $e_1$, we obtain 
 $2k^2 + 2x^2 + 2k -2x + 2y^2 + 1\leq \frac{1}{2}.$ 
However, the last inequality can not hold for any $x\in\mr, y\in\mr$, because $D_x=-4y^2-4k^2-4k<0$. 
Thus, the invertibility of  $M_{\overline{m},\Phi,\Psi}$ and $M_{m,\Psi,\Phi}$ can not be concluded from Proposition \ref{p3}.
\end{ex}

{\bf Proposition \ref{p4} applies, Proposition \ref{mpos} does not apply:}
\begin{ex}\label{8yes7no}
Let $\Phi=( e_1, e_1, e_2, e_3, e_4, \ldots )$, $m=(1,-1,1,1,1,1,\ldots)$, and $\Psi=((k+1)\phi_1, -(k+1)\phi_2, (k+1)\phi_3, (k+1)\phi_4,(k+1)\phi_5,\ldots)$, where $k\in(0,\frac{1}{2})$. By Example \ref{exdual}, the invertibility of $M_{m,\Phi,\Psi}$ and $M_{m,\Psi,\Phi}$ follows from Proposition \ref{p4}. Since $m$ is neither positive, nor negative, Proposition \ref{mpos} does not apply.
 \end{ex}

\v {\bf Proposition \ref{p3} applies, Proposition \ref{p4} does not apply:}

\begin{ex} \label{exdual2} Consider the frame $\Phi=( e_1, e_1, e_2, e_3, e_4, \ldots )$ with $A_\Phi^{opt}=1, B_\Phi^{opt}=2$ and the sequence $m\Psi=(e_2, e_1-e_2, e_2, e_3, e_4, \ldots)$, which is a dual frame of $\Phi$.
Clearly, Proposition \ref{p3} can be applied with $\Phi^d=m\Psi$ and arbitrary $\mu\in[0,\frac{1}{B_\Phi})$.

Now consider the sequence $m\Psi-\Phi$.
For the element $e_2$ we have 
$$\sum |\<e_2, m_n\psi_n - \phi_n\>|^2=|\<e_2,e_2-e_1\>|^2+|\<e_2,e_2\>|^2=2\|e_2\|^2.$$
Therefore, the sequence $m\Psi-\Phi$ can not be Bessel with bound 
$\mu<\frac{A_\Phi^2}{B_\Phi}$, because $\frac{(A^{opt}_\Phi)^2}{B^{opt}_\Phi}=\frac{1}{2}<2$.
 Thus, Proposition \ref{p4} can not be used. 
\end{ex}

\v {\bf Proposition \ref{p3} applies, Proposition \ref{mpos} does not apply:} 

\begin{ex} \label{9yes7no} Consider $\Phi=( e_1, e_1, e_2, e_3, e_4, \ldots )$, $m=(1,-1,1,1,1,1,1,\ldots)$ and $\Psi=(e_2, e_2-e_1, e_2, e_3, e_4, \ldots)$. By Example \ref{exdual2}, the invertibility of $M_{m,\Phi,\Psi}$ and $M_{m,\Psi,\Phi}$ follows from Proposition \ref{p3}.
Since $m$ is neither positive, nor negative, Proposition \ref{mpos} does not apply.
\end{ex}

{\bf Proposition \ref{mpos} applies, Proposition \ref{p4} does not apply:}

\begin{ex} \label{exnew}
Consider the frame $\Phi=(e_1,e_1,e_2,e_2,e_3,e_3,e_4,e_4,\ldots)$ with bounds $A_\Phi^{opt}=B_\Phi^{opt}=2$ and the frame $\Psi=(e_1,\frac{1}{2}\,e_1,e_2,e_2,e_3,e_3,e_4,e_4,\ldots)$. 
The sequence $\Psi-\Phi=(0,-\frac{1}{2}\,e_1,0,0,0,0,\ldots)$ is Bessel with $B_{\Psi-\Phi}^{opt}=\frac{1}{4}<2=\frac{(A_\Phi^{opt})^2}{B_\Phi^{opt}}$. Take $m=(4)$. Then $\frac{\sup |m_n|}{\inf |m_n|}=1<2\sqrt{2}=\frac{A_\Phi^{opt}}{\sqrt{B_{\Psi-\Phi}^{opt} B_\Phi^{opt}}}$. Now the invertibility of $M_{m,\Phi,\Psi}$ and $M_{m,\Psi,\Phi}$ can be concluded by Proposition \ref{mpos}. Since the sequence $m\Psi-\Phi=(3e_1,e_1,3e_2,3e_2,3e_3,3e_3,3e_4,3e_4,\ldots)$ is Bessel for $\h$ with bound $B_{m\Psi-\Phi}^{opt}=18>\frac{(A_\Phi^{opt})^2}{B_\Phi^{opt}}$, Proposition \ref{p4} can not be applied.

\end{ex}

{\bf Proposition \ref{mpos} applies, Proposition \ref{p3} does not apply:}

\begin{ex} \label{7yes9no}
Let $\Phi=( e_1, e_1, e_2, e_3, e_4, \ldots )$, $m=(1)$, and $\Psi=((k+1)\phi_n)$, where $k\in (0,\frac{1}{2})$. By Example \ref{exdual}, Proposition \ref{p3} can not be applied.
Again by Example \ref{exdual}, the sequence $\Psi-\Phi=(k\phi_n)$ is Bessel for $\h$ with optimal bound $B_{\Psi-\Phi}^{opt}=k^2B^{opt}_\Phi<\frac{(A^{opt}_\Phi)^2}{B^{opt}_\Phi}$.
Moreover, $\frac{\sup_n |m_n|}{\inf_n |m_n|}=1<\frac{1}{2k}=\frac{A_\Phi^{opt}}{\sqrt{B_{\Psi-\Phi}^{opt} B_\Phi^{opt}}}$.
 Thus, Proposition  \ref{mpos} implies the invertibility of $M_{m,\Phi,\Psi}$ and $M_{m,\Psi,\Phi}$.
\end{ex}

 {\bf Acknowledgments} 
The authors are thankful to the MULAC-partners for the fruitful discussions and valuable comments, in particular to H. Feichtinger, J.P. Antoine, F. Jaillet and D. Bayer, as well as J. White for proofreading. They are grateful for the hospitality of the 
Unit\'e de Physique Th\'eorique et de Physique Math\'ematique, Universit\'e Catholique de Louvain. 
The second author is also grateful for the hospitality of the Acoustics Research Institute and the support from the MULAC-project. She thanks 
 UACEG and the Department of Mathematics of UACEG for their support and openness  
  in order for the research on this paper to be done.

\bibliographystyle{plain}

\end{document}